\documentclass[a4paper]{IEEEtran}
\usepackage[latin1]{inputenc}
\usepackage[cmex10]{amsmath}
\usepackage{amsfonts}
\usepackage{amssymb}
\usepackage{graphicx}
\usepackage{float}
\usepackage{verbatim}
\usepackage{array}
\usepackage{multirow}
\usepackage{dcolumn}
\usepackage{color}
\usepackage[noadjust]{cite}
\usepackage{url}
\usepackage{balance}
\usepackage[usenames,dvipsnames]{xcolor}
\usepackage{algorithm}
\usepackage{algpseudocode}
\usepackage{breqn}
\usepackage{setspace}
\usepackage{amsmath}
\usepackage{diagbox}
\usepackage{adjustbox} 
\usepackage{makecell}
\usepackage{mathtools}
\usepackage{comment}

\usepackage{etoolbox}

\DeclareGraphicsExtensions{.eps}

\newcommand{\ignore}[1]{}

\definecolor{myblue}{gray}{0}

\begin{document}



\title{An Integrated Progressive Hedging and Benders Decomposition with Multiple Master Method to Solve the Brazilian Generation Expansion Problem}

\author{Alessandro~Soares,~Alexandre~Street,~Tiago~Andrade,~and~Joaquim~Dias~Garcia

\thanks{A. Soares and J. Dias Garcia are with PSR, Rio de Janeiro, Brazil, and with the Department of Electrical Engineering, Pontifical Catholic University of Rio de Janeiro, Brazil (\mbox{alessandro@psr-inc.com}, \mbox{joaquim@psr-inc.com})}

\thanks{A. Street is with the Department of Electrical Engineering, Pontifical
Catholic University of Rio de Janeiro, Brazil (\mbox{street@ele.puc-rio.br})}

\thanks{T. Andrade is with PSR, Rio de Janeiro, Brazil (\mbox{tiago@psr-inc.com})}
}
\maketitle


\begin{abstract}
This paper exploits the decomposition structure of the large-scale hydrothermal generation expansion planning problem with an integrated modified Benders Decomposition and Progressive Hedging approach. We consider detailed and realistic data from the Brazilian power system to represent hourly chronological constraints based on typical days per month and year. Also, we represent the multistage stochastic nature of the optimal hydrothermal operational policy through co-optimized linear decision rules for individual reservoirs. Therefore, we ensure investment decisions compatible with a nonanticipative (implementable) operational policy. To solve the large-scale optimization problem, we propose an improved Benders Decomposition method with multiple instances of the master problem, each of which strengthened by primal cuts and new Benders cuts generated by each master's trial solution. Additionally, our new approach allows using Progressive Hedging penalization terms for accelerating the convergence of the method. We show that our method is 60\% faster than the benchmark. Finally, the consideration of a nonanticipative operational policy can save 7.64\% of the total cost (16.18\% of the investment costs) and significantly improve spot price profiles.
\end{abstract}

\section*{Nomenclature}
\label{sec:nomenclature}

\subsection*{Sets and Indices}

\begin{IEEEdescription}[\IEEEusemathlabelsep\IEEEsetlabelwidth{$MMM$}]

\item[$\omega$]{ Index representing the operative scenario $\omega$ of the operation problem (subproblem).}
\item[$\Omega$]{Set containing the scenarios $\omega$.}
\item[$ \mathcal{X}$]{ Set {\color{myblue}constraining first stage decisions, i.e., investment decisions and linear decision rule coefficients}.}
\item[$\aleph_{t,d}$]{ Set of feasible hourly operating points modeling ramping and transmission capacity constraints at stage $t$ and typical day $d$.}
\item[$\tilde{V}_{t}$]{ Set of hydro operative constraints at stage $t$.}
\item[$\Lambda_{\omega}$]{ Feasibility set representing operational constraints.}
\item[$s$]{ Index representing an operative scenario selected to build a primal cut and identifying a master problem realization, its trial solution, and optimal value.}
\end{IEEEdescription}

\subsection*{Constants}

\begin{IEEEdescription}[\IEEEusemathlabelsep\IEEEsetlabelwidth{$MMM$}]
\item[$p_{\omega}$]{ Probability of scenario $\omega$.}
\item[$H_{t,d}$]{ Number of hours in typical day $d$, stage (month) $t$.}
\item[$c_v$]{ Vector of units' operational costs.}
\item[$I$]{ Vector of investment costs of generation units.}
\item[$R$]{ Incidence matrix allocating hydros' outflow to reservoirs according to rivers topology. {\color{myblue}Element $(i,j)$ is $-1$ if hydro $j$ is directly upstream of hydro $i$, diagonal elements are $1$, and all other elements are equal to zero}.}
\item[$a_{t}(\omega)$]{ Vector of hydro inflows at stage $t$ and scenario $\omega$.}
\item[$L$]{ {\color{myblue}Hydro selection matrix. Element $(i,j)$ is $1$ if generating unit $j$ is associated with hydro reservoir $i$, and $0$ otherwise}.}
\item[$P$]{ {\color{myblue}Diagonal coefficient matrix with the productivity of hydro generators}.}
\item[$G$]{ {\color{myblue}Diagonal matrix with maximum generation capacity of generating units}.}
\item[$A$]{ {\color{myblue}Generator-to-bus incidence matrix. Element $(i,j)$ is $1$ if generator $j$ is at bus $i$, and $0$ otherwise}.}
\item[$B$]{ {\color{myblue} Susceptance matrix of the DC power-flow approximation.}} 
\item[$D_{t,d,h}(\omega)$]{ Vector of demand per bus at stage $t$, typical day $d$, hour $h$ and scenario $\omega$.}
\item[$W,T$]{ Matrices modeling first- and second-stage coupling constraints.}
\item[$\mathcal{I}$]{Cost of first-stage variables.}
\item[$c$]{Cost for all operational variables.}
\item[$\rho$]{Progressive hedging static parameter.}
\item[$w_{s}^{k}$]{Progressive hedging dynamic parameter.}
\item[$\hat{x}_s^{k}$]{Candidate solution at scenario $s$ in the iteration $k$.}
\item[$\pi_s^{k}$]{Dual variable associated to the constraint that couples the investment and operation problem, at scenario $s$ and iteration $k$.}
\item[$\bar{x}^{k}$]{Average of the first-stage decision from iteration $k-1$.}
\end{IEEEdescription}

\subsection*{Decision Variables}

\begin{IEEEdescription}[\IEEEusemathlabelsep\IEEEsetlabelwidth{$MMM$}]
\item[$g_{t,d,h}(\omega)$]{ Generation vector at stage $t$, typical day $d$, hour $h$, and scenario $\omega$.}
\item[$x^{INV}$]{ Binary vector of investment decisions (existing generators are modeled as fixed entries equal to 1).}
\item[$v_{t}(\omega)$]{ Vector of reservoirs storage at the end of stage $t$ and scenario $\omega$.}
\item[$u_{t}(\omega)$]{ Vector of amount of water used to generate electricity during stage $t$ and scenario $\omega$.}
\item[$s_{t}(\omega)$]{ Vector of water spilled during stage $t$ and scenario $\omega$.}
\item[$v_{0}(\omega)$]{ Vector of initial storage condition at stage $t$ and scenario $\omega$.}
\item[$\theta_{t,d,h}(\omega)$]{Vector of buses' phase angle at stage $t$, typical day $d$, hour $h$ and scenario $\omega$.}
\item[$x_{t}^{LDR}$]{ Linear decision rule vector of angular coefficients for stage $t$.}
\item[$x_{t,0}^{LDR}$]{ Linear decision rule vector of linear coefficients per hydro unit for stage $t$.}
\item[$Q(\cdot)$]{ Real operational cost function.}
\item[$x$]{ Vector of first-stage decision variables comprising the binary investment variables and linear decision rule coefficients.}
\item[$y(\omega)$]{ Vector of decision variable representing all operative variables at scenario $\omega$.}
\item[$\pi_{\omega}$]{ Vector of dual variables associated with operation-investment coupling constraints at scenario $\omega$, used to build the Benders cuts.}
\item[$\alpha_{\omega}$]{ Approximation of the real operational cost at scenario $\omega$.}
\end{IEEEdescription}

\begin{table*}[h]
\caption{Proposed approach compared to literature}
\label{tab:literature-compare}
\resizebox{\textwidth}{!}{%
\begin{tabular}{ccccccc}
\hline
\textbf{Approach}  & 
\textbf{\begin{tabular}[c]{@{}c@{}}Representation of \\ Uncertainties\end{tabular}} &
\textbf{\begin{tabular}[c]{@{}c@{}}Hourly \\ constraints\end{tabular}} & \textbf{\begin{tabular}[c]{@{}c@{}} Operational\\ policy\end{tabular}}                     & \textbf{\begin{tabular}[c]{@{}c@{}}Binary investment \\ decision\end{tabular}} & \textbf{\begin{tabular}[c]{@{}c@{}}Co-optimization of \\ energy and reserves\end{tabular}} & \textbf{\begin{tabular}[c]{@{}c@{}}Decomposition \\ Technique\end{tabular}} \\ \hline
Koltsaklis, N. E. et al. (2015) [17] \ignore{\cite{koltsaklis2015multi}}  & -                                                                                                                  & \checkmark                                                          & -                                                                                                                      & \checkmark                                                      & \checkmark                                                                  & -                                \\
Pina, A. et al. (2013) [20] \ignore{\cite{pina2013high}}   & -                                                                                                                  & \checkmark                                                          & -                                                                                                                      & \checkmark                                                      & -                                                                                          & Feasibility cuts                 \\
Li, J. et al. (2018) [21] \ignore{\cite{li2018robust} }                       & \checkmark $^A$ & \checkmark                                                          & -                                                                                                                      & \checkmark                                                      & \checkmark                                                                  & Column-and-Constraint            \\
Liu, Y. et al. (2018) [22] \ignore{\cite{liu2018multistage}}& \checkmark                                                                                          & \checkmark                                                          & -                                                                                                                      & -                                                                              & -                                                                                          & Progressive Hedging              \\
Thome, F. et al. (2019)  [23] \ignore{\cite{thome2019stochastic}}                       & \checkmark                                                                                          & -                                                                                  & \checkmark                                                                                              & \checkmark                                                      & \checkmark                                                                  & Benders Decomposition            \\
Proposed approach                               & \checkmark                                                                                          & \checkmark                                                          & \checkmark $^B$ & \checkmark                                                      & \checkmark                                                                  & BDPH                             \\ \hline
\end{tabular}}\\[.05cm]
\footnotesize{$^A$ The uncertainties are represented through robust optimization methods}\\[.05cm]
\footnotesize{$^B$ The hydro operation policy is represented through linear decision rules}
\end{table*}

\section{Introduction}


\IEEEPARstart{G}{eneration} expansion planning (GEP) models aim to minimize the total cost of investment and operation through a long-term horizon. They bring relevant insights for market agents and also provide planners and regulators with valuable information about long-term equilibrium generation portfolios under the absence of market power abuse \cite{munoz2018economic,munoz2017does}. This is especially relevant in hydrothermal power systems, where important regulatory and economic metrics needed to induce efficient generation expansion rely on the assessment of the opportunity cost of water under a long-term investment and operational equilibrium. For instance, the current regulatory guidelines of the Brazilian power system, relies on the definition of an indicative generation expansion plan to rank generators' offers in new generation auctions. This ranking process is based on a (\$/MWh)--cost-benefit metric calculated with long-run (equilibrium) spot-price scenarios simulated for the indicative expansion plan. {\color{myblue}The indicative plan is obtained by a simplified GEP model that disregards the effect of short-term uncertainties and constraints as well as the co-optimization of investment decisions and the nonanticipative scheduling of reservoirs.} For the interested reader, we refer to \cite{Bezerra2006} for a short paper on the assessment of the cost-benefit metric in Brazilian new generation auctions, to \cite{maurer2011electricity} for a more comprehensive reading about related energy auctions and mechanism used in Latin America, and to \cite{EPE2020} for the official 2020 expansion-plan report delivered by the Brazilian system planner. 

GEP models should take into account the main long- and short-term characteristics of the system and uncertainties to enhance the description of {\color{myblue}operational} opportunity costs {\color{myblue}and thereby the quality} of first-stage decisions \cite{Moreira2016,brigatto2017assessing}. Furthermore, this class of models is usually non-convex due to the necessity of representing integer investment and operative decisions. Conventional approaches to solve these large-scale non-convex stochastic programs are: non-linear programming (NLP); mixed-integer linear programming (MILP); and decomposition techniques (such as Benders decomposition). All of these approaches may be used together with approximations and assumptions to make the problem computationally tractable \cite{ringkjob2018review}. {\color{myblue} In this context, the search for more efficient decomposition algorithms is key for the development of more realistic expansion planning models.}

A wide range of applications on GEP are found in the literature \cite{koltsaklis2018state,oree2017generation}, each of which considering different aspects and system characteristics. In this work, however, we focus on the specific challenges of hydrothermal power systems \cite{maluenda2018expansion,costa2021reliability}. In this setting, the main challenge is to consider, within the expansion problem, an integrated and computationally efficient nonanticipative water-value assessment, {\color{myblue}a rarely explored subject in the related technical literature as will be further depicted}. Therefore, hydrothermal-based GEP largely relies on the co-optimization of investment decisions and long-term multistage stochastic dispatch policies. Also, there are several other applications with a similar structure, such as the maintenance optimization problem \cite{helseth2018optimal}, for which the algorithms and ideas exploited in this work are also valid.

\subsection{Energy resource planning models}

In energy resource planning literature, \cite{pineda2018chronological} solves the investment and operation problems simultaneously. To reduce the computational burden and avoid decomposition, the authors adopt a linear relaxation of investment decisions and propose {\color{myblue}a clustering method}
to reduce the number of hours yet keeping a chronological representation of externalities. A clustering algorithm is also proposed in \cite{liu2018hierarchical}, where the author shows the strategy's effectiveness comparing clustered problems with the unclustered version. In \cite{poncelet2016selecting}, the authors propose a novel optimization problem to minimize the approximation errors of the typical days. Time clustering schemes are commonly used in the literature as typical days or weeks. For instance, \cite{koltsaklis2015multi} uses typical days per month and year, and considers a deterministic model for long-term planning with a detailed representation of short-term constraints. In this setting, the problem can be solved without a decomposition approach. \cite{maluenda2018expansion} proposes a hydrothermal-based GEP, using typical days to represent hourly constraints and scenarios to represent hydro inflow uncertainties. The authors use the progressive hedging (PH) technique to decompose the problem scenario-wise. Several other recent planning models in literature adopt simplifications to make the investment and operational problems computationally tractable \cite{rashidaee2018linear,zhang2018mixed}. 

Notwithstanding, decomposition methods are largely applied in this subject. For instance, \cite{pina2013high} proposes a decomposition approach where the master consists of a deterministic investment problem and the subproblem is a detailed short-term operational problem that produces feasibility cuts. \cite{li2018robust} addresses uncertainties in the net load by a robust optimization approach, considering uncertainty in the hourly ramping, but without other detailed short-term constraints. \cite{liu2018multistage} proposes a multiscale multistage stochastic model, addressing short-term constraints and uncertainties, and decomposes the problem with PH, where the problem is convex (with linear investment and commitment decisions), which guarantees optimal solutions but may not be as fast as BD. Table \ref{tab:literature-compare} summarizes the comparison between the proposed approach and the related energy resource planning literature \cite{koltsaklis2015multi,pina2013high,li2018robust,liu2018multistage,thome2019stochastic}. In this table, symbols "\checkmark" and "-" indicate whether a particular aspect is considered or not.

\subsection{Decomposition structures}

Deterministic planning models can be solved in a reasonable amount of time, even considering the co-optimization of investment and operational decisions. Uncertainty representation drastically increases the size of the problem, leading to intractability issues, especially when the number of scenarios is large. This is the case for most real problems. Notwithstanding, these issues can be especially worsened in the presence of time-coupling constraints requiring a multistage model to characterize the opportunity costs of operational resources such as water. Hence, decomposition approaches such as PH \cite{rockafellar1991scenarios} and Benders Decomposition (BD) \cite{benders1962partitioning} are frequently used. PH algorithms guarantee convergence to the optimal solution when the problem is convex. However, since real expansion problems have binary investment variables, PH is usually used as a heuristic to obtain solutions \cite{munoz2015scalable}.
{\color{myblue} BD techniques guarantee optimal solutions when first-stage decisions are mixed integer and the recourse function is convex.}
The BD approach was first proposed in the context of GEP by Campodonico et al. \cite{campodonico2003expansion} and is used to solve many real problems \cite{rebennack2014generation,latorre2019stochastic,oliveira2007value} since the optimal solution is guaranteed in a reasonable amount of time.

The Progressive Hedging (PH) was proposed by Rockafellar and Wets \cite{rockafellar1991scenarios}. It is a decomposition approach based on the augmented Lagrangian Relaxation. Watson and Woodruff \cite{watson2011progressive} describe practical details of PH to a class of large scale stochastic mixed-integer resource allocation problems. 
Guo et al. \cite{guo2015integration} used the PH to speed up a dual-decomposition-based algorithm, whereas \cite{gade2016obtaining} proposed a PH lower bound for two-stage and multi-stage stochastic mixed-integer programs. 

Crainic et al. \cite{crainic2016partial} proposed the partial BD, where they add a subset of constraints and variables to the master problem. The algorithm we present in this work improves this idea by considering multiple copies of the master problem, each of which accounting for primal cuts generated based on scenario information and new Benders cuts obtained from each master's trial solution. Due to the multiplicity of trial solutions generated by the information of different scenarios, we use Progressive Hedging penalization terms to increase the consensus among trial solutions and accelerate our method. 

\vspace{+0.2cm}
\subsection{Nonanticipative hydrothermal dispatch}

The solution of a long-term GEP problem applied to a hydrothermal power system requires the consideration of a multistage reservoir operational policy. The main reason for that is to avoid the threat of optimistically biasing the water opportunity-cost assessments based on an anticipative operational model \cite{pereira1991multi,street2020assessing}. In other words, we need to consider in our GEP a decision rule (see \cite{shapiro2014lectures}) that is as close as possible to an implementable (nonanticipative) policy to avoid under investments due to artificially reduced operational costs (based on optimistic anticipative policies) that will not be achievable in practice. In this context, the customary two-stage approximation -- in which given the investment decisions, the system operation follows with perfect (anticipative) information of the uncertainty realizations (i.e., per scenario) -- is not valid as we demonstrate in our case study. 

A nonanticipative operational policy is a rule defining decision variables of a given period $t$ based on previously reveled information, i.e., without assuming access to the information of uncertainty factors after $t$. The linear decision rule (LDR) methodology defines an implementable nonanticipative policy based on an optimized linear combination of functions applied to the previously revealed uncertainty scenario. \cite{revelle1969linear} first proposed the LDR for reservoir management, and this approach is gaining more and more attention in the literature (see \cite{chen2008linear,egging2016linear,braaten2016linear} and \cite{bodur2018two}). {\color{myblue} Notwithstanding, to the best of the authors knowledge, the consideration of LDR to address the relevant piece of realism of considering nonanticipative multistage stochastic dispatch policies within a hydrothermal GEP model has not being addressed before. Thus, in} this work, we propose a new stochastic hydrothermal GEP with multistage (nonanticipative) dispatch policy based on LDRs.



\subsection{Contributions and work organization}
The main contributions of this work are threefold:
\begin{itemize}
    \item {\color{myblue}A new hydrothermal GEP model considering nonanticipative multistage stochastic dispatch policies through linear decision rules.} Out-of-sample tests based on real data demonstrate that the consideration of {\color{myblue}more realistic nonanticipative dispatch policies} has significant impacts on first-stage investment decisions and subsequent operation costs {\color{myblue}and spot-price profiles}.
    \item An improved Benders Decomposition with multiple master problems (BDMM) method, each of which strengthened by primal cuts based on scenario information and new Benders cuts generated by each master's trial solution.
    \item Leveraging the diversity of trial solutions our BDMM method provides, we combine the proposed BDMM with Progressive Hedging penalization terms to find a consensus among the multiple solutions and accelerate the proposed method. Therefore, we propose a novel accelerated Benders Decomposition with Multiple Masters (a-BDMM) method based on the Progressive Hedging (PH) consensus idea.
\end{itemize}

The remainder of this paper is organized as follows. In Section \ref{sec:int}, the GEP problem is introduced and formulated. The decomposition algorithm is presented as the solution strategy to solve the problem in Section \ref{sec:strategy}. Section \ref{sec:results} provides numerical results illustrating the performance of the proposed algorithm and an analysis of the anticipative policy in the GEP formulation. Finally, in Section \ref{sec:conclusions} final remarks are drawn.


\section{The generation expansion planning model}
\label{sec:int}

This section introduces the GEP problem formulation as a MILP optimization problem, which will be referred to as the Deterministic Equivalent (DE)  problem. We assume a discrete and finite sample space $\Omega = \{1,..,\omega,...\}$, in which each scenario $\omega$ is assumed to have a known conditional probability $p_{\omega}$. For the operational variables, we use a LDR to consider a monthly nonanticipative multistage operational policy under uncertainty of inflows \cite{shapiro2005complexity, bodur2018two} for the reservoirs.

The system operation constraints and costs are computed within an hourly resolution based on monthly hydro generation targets dictated by the LDR. Thus, based on monthly inflow scenarios and subsequent hydro generation amounts given by the LDR, the uncertainties of intermittent renewable sources are used to define the operation of typical (representative) days within an hourly granularity. In this sense, we approximate the daily operation within each stage (months) by weekdays, weekends, and holidays multiplied by their number of hours per month. Hourly scenarios for renewables are conditionally generated based on the inflow scenarios to present correlations.

Mathematically, for each month (stage) $t\in \mathcal{T}$ and scenario $\omega\in\Omega$ we have: i) strategic stagewise decisions, such as $u_t(\omega),s_t(\omega),v_t(\omega)$, defining the total amount of water used to generate electricity and spilled during the stage, and the storage level at the end of the stage; ii) short-term operational decisions, such as $g_{t,d,h}(\omega), \theta_{t,d,h}(\omega)$, defining the hourly generation and the bus angles for each typical day $d\in \mathcal{D}$ and hour $h\in \mathcal{H}$ of the stage $t$; iii) constraints linking strategic and short-term decisions, $Pu_{t}(\omega) = L\sum_{d,h}H_{t,d}g_{t,d,h}(\omega)$; and iv) $x^{INV}$ and $(x_{t,0}^{LDR},x_{t}^{LDR})$ representing the investment decisions and LDR coefficients (first-stage decision vectors). The proposed model is detailed as follows in expressions \eqref{mod:init}-\eqref{mod:invcstr}:
\begin{align}
   &\min \ \  I^{\top}x^{INV} + \sum\limits_{t,d,h,\omega}p_{\omega}c_v^{\top}H_{t,d}g_{t,d,h}(\omega)\label{mod:init}\\
               &s.t.   \notag  \\
                       &g_{t,d,h}(\omega) \leq Gx^{INV} \ \ \forall t,d,h,\omega \label{mod:invcouple}\\
                       &u_{t}(\omega) = x_{t,0}^{LDR} +  x_{t}^{LDR} a_{t}(\omega) \ \ \forall t,d,h,\omega \label{mod:affine}\\
                       &v_{t}(\omega) = 
                 v_{t-1}(\omega) + a_{t}(\omega) - R(u_{t}(\omega)+s_{t}(\omega)) \ \ \forall t,\omega \label{mod:hydbal}\\
                       &v_{T}(\omega) \ge v_{0}(\omega) \ \ \forall t,\omega\label{mod:vivf}\\
                       &Pu_{t}(\omega) = L\sum\limits_{d,h}H_{t,d}g_{t,d,h}(\omega) \ \ \forall t,\omega \label{mod:hidcouple}\\
                       &Ag_{t,d,h}(\omega) + B\theta_{t,d,h}(\omega)= D_{t,d,h}(\omega) \ \ \forall t,d,h,\omega \label{mod:loadbalance}\\
                       &\{ g_{t,d,h}(\omega), \theta_{t,d,h}(\omega) \}_h \in \aleph_{t,d} \ \ \forall t,d,\omega\label{mod:opecstr}\\
                       &\{ v_{t}(\omega), u_{t}(\omega),s_{t}(\omega) \}\in \tilde{V}_{t} \ \ \forall t,\omega. \label{mod:hydro} \\
                       &\{x^{INV}, x_{t,0}^{LDR},  x_{t}^{LDR}\} \in \mathcal{X}  \label{mod:invcstr}
\end{align}




{\color{myblue}
We present a compact vector formulation using matrices defined in the nomenclature section.} We omit the sets in which indexes range for the sake of conciseness\footnote{{\color{myblue} We refer the interested reader to \cite{soares2019optimal} for an extensive and detailed formulation of the operational model that is similar to the one used in this paper.}}. The objective function \eqref{mod:init} has two parts, the investment cost and the present value of the expected operational costs.
Constraint \eqref{mod:invcouple} represents the relation between investment and operational decisions, {\color{myblue}i.e., a generating unit is only available to produce energy if the investment cost was paid.}
Constraint \eqref{mod:affine} represents the linear decision rule. In this expression, $x_{t,0}^{LDR}$ and $x_{t}^{LDR}$ are decision variables representing the vector of linear and matrix of angular coefficients defining the LDR, respectively. {\color{myblue}Hence, turbine outflows must be affine functions of the inflows.}
{\color{myblue}Constraint \eqref{mod:hydbal} describes for each reservoir the water mass balance equation, where the final volume of each reservoir at the end of stage $t$ is equal to its initial volume, plus the incoming inflow, minus the total net water discharged downstream. Note that the matrix $R$ accounts for the hydro cascades as per its definition in the nomenclature section.}
Constraint \eqref{mod:vivf} represents the hydro operational strategy used to prevent the end-of-horizon effect. This constraint aims to obligate the model to use only the water that arrives along the years in the study horizon, {\color{myblue}representing a cyclic and sustainable usage of the water resources.}
Constraint \eqref{mod:hidcouple} refers to the hourly hydro generation modeling, where an average productivity approximation is considered for each hydro as widely adopted in long-term studies (see relevant publications in the last three years \cite{rosemberg2021assessing,lohndorf2019modeling,street2020assessing,debia2021strategic}).
{\color{myblue}Constraint \eqref{mod:loadbalance} refers to the load balance, in which the sum of generation in each bus plus net energy transfer must match the demand. Note that power flows are represented by angle differences multiplied by line sucesptance following the widely used DC power-flow approximation for the AC network model \cite{rosemberg2021assessing}.}
Constraint \eqref{mod:opecstr} represents other operational constraints, such as ramping constraints and limits on angle {\color{myblue}differences} modeling transmission lines maximum flow capacity, {\color{myblue}suppressed to keep the formulation compact.} {\color{myblue}Constraint \eqref{mod:hydro} represents physical bounds on the hydro plant variables}, such as minimum and maximum values for water storage, maximum generation and spillage.
{\color{myblue}Finally, expression \eqref{mod:invcstr} represents generic constraints on the first-stage variables (investment variables and coefficients of the linear decisions rules). For the case of investment variables, \eqref{mod:invcstr} allows considering target capacity, energy policies, and other constraints on investment decisions. Additionally, it {\color{myblue} also allows for the consideration of constraints on} the LDR coefficients{\color{myblue}. Notwithstanding, in this work we only considered the binary nature of investment decisions}.}

%
\section{Solution Strategy}
\label{sec:strategy}
Problem \eqref{mod:init}-\eqref{mod:invcstr} provides optimal investment decisions considering a multistage stochastic operational policy. However, because we are using a LDR parametrization, the problem can be formulated as a large-scale multiperiod two-stage stochastic optimization problem, where the LDR is estimated as part of first-stage variables. To efficiently solve this problem, in the next sections, we present a variant of the Benders decomposition approach. 
For the sake of simplicity and didactic purposes, we re-write problem \eqref{mod:init}-\eqref{mod:invcstr} in a compact formulation as follows:
\begin{align}
   \min \ \ \ \ \ &\mathcal{I}^{\top}x + \sum\limits_{\omega\in\Omega}p_{\omega}c^{\top}y_\omega\label{eq:bigstart}\\
               s.t.    &\ \ \ \ Wy_\omega - Tx \leq h_\omega \ \ \forall \omega\in\Omega\ \label{eq:coupling}\\
                       &\ \ \ \ x\in\mathcal{X}, \ y_\omega \in \Lambda_{\omega} \ \ \forall \omega\in\Omega, \label{eq:bigend}
\end{align}

\noindent where, $y_\omega$ is the vector comprising the operational variables ($g_{t,d,h}(\omega)$, $\theta_{t,d,h}(\omega)$,  $v_{t}(\omega)$,  $u_{t}(\omega)$, $s_{t}(\omega)$); $\Lambda_w$ is a set containing the feasibility constraints \eqref{mod:hydbal}-\eqref{mod:hydro}; $x$ is the vector comprising the first-stage variables ($x^{INV}$, $x_{t,0}^{LDR}$, $x_{t}^{LDR}$); and constraint \eqref{eq:coupling} couples the first- and second-stage variables corresponding to \eqref{mod:invcouple}-\eqref{mod:affine}.

\vspace{-0.2cm}
\subsection{Traditional Benders decomposition}\label{sec:tbd}
Now we present the traditional Benders decomposition (TBD) applied to solve problems \eqref{eq:bigstart}-\eqref{eq:bigend}. We decompose the problem into master and subproblem problems. The master selects the first-stage variables (vector $x$ comprising investment and LDR coefficients), while the subproblem evaluates the recourse function, $Q(x)$, by solving the operative problem given $x$. 

So, we start defining the recourse function (the expected cost of the second-stage) for a given point $x$ as follows:
\begin{equation}\label{eq:recourse}
    Q(x) = \sum\limits_{\omega\in\Omega}p_{\omega}q_{\omega}(x).
\end{equation}
The evaluation of $Q(x)$ can be decomposed per scenario $\omega$ and solved in parallel. For each scenario, $q_\omega(x)$ represents the second-stage cost and can be calculated as follows:
\begin{align}
   q_\omega(x) = \min & \ c^{\top}y_\omega \label{eq:opestart}\\
               s.t.     & \  Wy_\omega \leq Tx + h_\omega \ :\xrightarrow[]{\text{{\small dual}}}\ \pi_{\omega} \label{eq: coupling constraint} & \ \\
                        &\ \ \ \ y_\omega \in \Lambda_{\omega}, \ &  \label{eq:opeend}
\end{align}
where {\color{myblue}$\pi_{\omega}$ is the dual of constraint~\eqref{eq: coupling constraint}, and $\pi_{\omega}^{\top}T$ is a subgradient of $q_\omega$ with respect to $x$}. Then, for a given iteration $k$ of the algorithm, we run the master problem to obtain a new trial solution, $\hat{x}^{k}$, and a lower bound, $LB_{\text{TBD}}$. The master is a relaxation of problem \eqref{mod:init}-\eqref{mod:invcstr} because the recourse function is approximated from below by supporting planes. These planes are also called Benders cuts and are obtained in previous iterations of the method. Using the multi-cut method \cite{rahmaniani2017benders}, the master problem of a TBD returns a new trial solution and a lower bound as follows: 
\begin{align}
   z^k, \hat{x}^{k} &\leftarrow \min \ \ \mathcal{I}^{\top}x + \sum\limits_{\omega\in\Omega}p_{\omega}\alpha_{\omega} \label{eq:masterinvstart}\\
               s.t.\ & \alpha_{\omega} \geq q_{\omega}(\hat{x}^{j}) + (\pi_\omega^{j})^{\top}T(x-\hat{x}^{j})\ \nonumber \\
               & \hspace{3cm} \forall \omega\in\Omega, j\in[k-1]\\
                 \ \ &\alpha_{\omega} \geq 0\ \forall \omega\in\Omega, \ x\in\mathcal{X},\label{eq:masterinvend}
\end{align}
\noindent where $\alpha_{\omega}$ represents the best approximation of the epigraph of $q_{\omega}$ until iteration $k$. Furthermore, hereinafter, we adopt the notation in which $[k-1]=\{1,...,k-1\}$ and $[0]=\emptyset$. Then, by solving \eqref{eq:opestart}--\eqref{eq:opeend} for the newly obtained trial solution $\hat{x}^{k}$, a new Benders cut can be generated to feed the next iteration master problem. Additionally, a lower and upper bound can be assessed to check the optimality GAP of the current solution as follows:
\begin{align}
&LB_{\text{TBD}} = z^k\\
&UB_{\text{TBD}} = \mathcal{I}^{\top}\hat{x}^{k} + Q(\hat{x}^{k}).
\end{align}
If the $GAP=UB_{\text{TBD}}-LB_{\text{TBD}}\le\epsilon$, then the algorithm stops and $x^k$ is returned as the optimal solution. Otherwise, the $k$ is incremented, and the master problem is called once again.

\vspace{-0.2cm}
\subsection{The Benders Decomposition with multiple master problems}
In this section, we present our proposed BDMM method. Firstly, we make $S=|\mathcal{S}|$ copies of the master problem \eqref{eq:masterinvstart}-\eqref{eq:masterinvend}, i.e., in each iteration of our BDMM method we define one master problem for each $s\in\mathcal{S}$. Secondly, each master problem, now indexed by $s$, instead of considering only Benders cuts, also considers the second-stage primal constraints associated with a given scenario $s$. It is worth highlighting that $\mathcal{S}$ could be generated based on different clustering strategies. Hence, we will develop our method for a general set $\mathcal{S}$. However, in this paper, we will use $\mathcal{S}=\Omega$, {\color{myblue}which may lead us to interpret the multiple masters as a stochastic master problem}. Then, we end up with $S=|\Omega|$ master problems, each of which differing from each other by a (stochastic) primal cut related to a given scenario. Furthermore, due to the multiplicity of master problems, multiple ($S$) trial solutions are also generated. Therefore, for each one of the $S$ newly generated points $\{x^k_{s}\}_{s\in\mathcal{S}}$, the multi-cut approach generates $|\Omega|$ new cuts, each of which approximating one function in $\{q_\omega(\cdot)\}_{\omega\in\Omega}$. Consequently, in each master problem $s$ of a given iteration $k$, a total of $S\cdot|\Omega|$ Benders cuts are considered for each previous iterations. Thus, the $s$--master problem and the associated lower bound and trial solution are defined as follows:
\vspace{-0.4cm}

\begin{align}
       &z^k_s, \hat{x}^k_s \leftarrow \min   \ \mathcal{I}^{\top}x + \sum\limits_{\omega}p_{\omega}\alpha_{\omega} \label{eq:masterimprovedstart} \\
       &s.t.\ \alpha_{\omega} \geq q_{\omega}(\hat{x}_{s'}^{j}) + (\pi_{\omega,s'}^{j}){\color{myblue}^{\top}T}(x-\hat{x}_{s'}^{j})\nonumber\\ 
           & \hspace{3cm}\forall \omega \in\Omega, s'\in S, j\in[k-1] \label{eq:masterimprovedendTBC}\\
           & \ \alpha_{s} \geq c^{\top}y_s \ \label{eq:masterimprovedendcut}\\
        & \ Wy_s - Tx \leq h_s \label{eq:masterimprovedendcut2} \\
           & \ y_s \in \Lambda_s \label{eq:masterimprovedendcut3} \\  
           & \ x\in\mathcal{X}. \label{eq:masterimprovedend}
    \end{align}

Note that traditional Benders cuts are built based on local-dual information of the recourse problem, thereby providing the master problem with loose linear approximations of the recourse function. The primal cut defined by \eqref{eq:masterimprovedendcut}--\eqref{eq:masterimprovedendcut3}, on the other hand, provides a much richer polyhedral information about the second stage to the master problem. This improvement proposed in this work is inspired by the success of column-and-constraint-generation algorithms applied to robust optimization, where few primal cuts are actually needed to support the optimal decisions (see \cite{street2013energy}). Furthermore, the number of Benders cuts in \eqref{eq:masterimprovedendTBC} is $S$ times greater than in the TBD method, which significantly improves the description of the recourse function.


After solving the $S$ instances of the master problem, we have $\{z^k_s, \hat{x}^k_s\}_{s\in\mathcal{S}}$. Because all values in $\{z^k_s\}_{s\in\mathcal{S}}$ are valid lower bounds for the problem, an improved Benders lower bound can be calculated based on the maximum among all values, i.e.,
\begin{align}
&LB_{\text{BDMM}} = \max\limits_{s\in\mathcal{S}} \{z^k_s\}.\label{BDMMLB}
\end{align}
A similar approach can be used to improve the upper bound. By evaluating the recourse function on each $\hat{x}^k_s$, we get $S$ new candidates for upper bounds, $\{\mathcal{I}^{\top}\hat{x}^k_s + Q(\hat{x}^k_s)\}_{s\in\mathcal{S}}$. Thus, we can select the lowest upper bound, i.e., we define
\begin{align}
&UB_{\text{BDMM}} = \min\limits_{s\in\mathcal{S}}\{ \mathcal{I}^{\top}\hat{x}^k_s + Q(\hat{x}^k_s)\}\label{BDMMUB},
\end{align}
\noindent and store the solution associated with the best upper bound, $\hat{x}^k(s^*)$, as the best trial solution at iteration $k$. A comparison with the best solution found so far is also advisable to keep the global best solution at hand. 

Finally, note that a lower and upper bound comparison between the TBD and the proposed BDMM is not directly possible because the two methods should follow different paths. However, it is clear that the BDMM provides tighter approximations in every master problem since the first iteration. The improvement in lower bounds comes at the cost of a higher computational effort. The tradeoff between improving the lower bound assessment and the additional computational effort will be depicted in our Case Study Section. Leveraging the diversity of trial solutions our BDMM method provides, in the next section, we propose a novel acceleration scheme.


\vspace{-0.3cm}
\subsection{Accelerating convergence of the Benders Decomposition with multiple masters}
\vspace{-0.2cm}
Now we present a new consensus scheme based on the PH method to accelerate our BDMM. To do that, we add the PH penalty terms in the master problem formulation \eqref{eq:masterimprovedstart}-\eqref{eq:masterimprovedend}. So, we rewrite the master problem \eqref{eq:masterimprovedstart}-\eqref{eq:masterimprovedend} as problem \eqref{eq:master2start}-\eqref{eq:master2end}, adding the terms related to the augmented Lagrangian relaxation following the PH approach (see \cite{gade2016obtaining}). 

\begin{align}
   &\hat{x}^k_s \leftarrow \min \ \ \mathcal{I}^{\top}x + \sum\limits_{\omega}p_{\omega}\alpha_{\omega} +  \nonumber\\
   &\hspace{2.5cm}\frac{\rho}{2}\|x-\bar{x}^{k}\|^2 + w_{s}^{k\top}(x-\bar{x}_k) \label{eq:master2start} \\
   \ \ &s.t.
   \ \ \text{Constraints \eqref{eq:masterimprovedendTBC}--\eqref{eq:masterimprovedend}.\label{eq:master2end}}
\end{align}
\noindent
Then, after solving the $S$ master problems, following the BDMM approach, $w_{s}^{k}$ is updated following the sub-gradient method, i.e.,
\begin{equation}
\label{eq:wupdate}
w_{s}^{k+1} = w_{s}^{k} + \rho(\hat{x}^{k}_s - \bar{x}^{k})
\end{equation}
\noindent
where $\hat{x}^{k}_s$ is the solution of the problem \eqref{eq:master2start}-\eqref{eq:master2end} from the iteration $k$ and scenario $s$, and $\bar{x}^{k}$ is the average of all the $S$ trial solutions obtained with  \eqref{eq:master2start}-\eqref{eq:master2end} at iteration $k$.

To obtain a lower bound, however, we have to solve the following modified version of the problem, which considers only the simple Lagrangian relaxation (without the quadratic terms):  
%
\begin{align}
   \zeta^k_s = & \min \ \ \ \mathcal{I}^{\top}x + \sum\limits_{\omega}p_{\omega}\alpha_{\omega} + w_{s}^{k\top}(x-\bar{x}_k) \label{eq:masterLBstart} \\
    \ \ &s.t. \ \ \text{Constraints \eqref{eq:masterimprovedendTBC}--\eqref{eq:masterimprovedend}}. \label{eq:masterLBend}
\end{align}
Thus, the lower bound can be calculated as follows:  \\
\textbf{Theorem 1.} \emph{Problem \eqref{eq:bigstart}-\eqref{eq:bigend} admits the following lower bound:}
\begin{align}
&LB_{\text{a-BDMM}} = \sum\limits_{s=1}^{S} p_s\zeta^k_s.
\end{align}
Although slightly different, the proof to \textbf{Theorem 1} goes very much like that provided in \cite{gade2016obtaining}. Therefore, for the sake of conciseness, we omit the proof here. Finally, the upper bound remains unchanged and follows expression \eqref{BDMMUB}. The proposed a-BDMM algorithm is summarized as follows: 
\label{sec:algorithm}
\begin{algorithm}[h]
  \small
\caption{The proposed a-BDMM method}
\label{alg:alg-a-BDMM}
\begin{algorithmic}[1]
    \State \textbf{Initialization ($\rho \leftarrow \textbf{input}$)} 
    \State \ \ $k\leftarrow 0, 
    GAP^k \leftarrow +\infty, w_{s}^k\leftarrow0 \ \ \forall s\in \mathcal{S}$ 
    \While{$GAP^k>\epsilon$}
        \State $k\leftarrow k+1$
        \For{each $s\in\mathcal{S}$ (computed in parallel) }
            \State Solve master problem \eqref{eq:master2start}--\eqref{eq:master2end} and store $x^k_s$
            \State \textbf{for} each $\omega\in\Omega$: compute $q_w(x^k_s)$ and store $\pi^{k}_{\omega,s}$
            \State Solve problem \eqref{eq:masterLBstart}--\eqref{eq:masterLBend} and store $\zeta^k_s$
        \EndFor
        \State \textbf{Compute:} 
            \State \ \ \ \ $\bar{x}^{k} \leftarrow \sum\limits_{s\in\mathcal{S}}p_{s}\hat{x}^{k}_s$
            \State \ \ \ \ $LB^k_{\text{a-BDMM}} \leftarrow \sum\limits_{s=1}^{S} p_s\zeta^k_s$
            \State \ \ \ \ $UB^k_{\text{a-BDMM}} \leftarrow \min\limits_{s\in\mathcal{S}}\{ \mathcal{I}^{\top}\hat{x}^k_s + Q(\hat{x}^k_s)\}$\vspace*{0.1cm}
            \State \ \ \ \ $GAP^k \leftarrow UB^k_{\text{a-BDMM}}- LB^k_{\text{a-BDMM}}$\vspace*{0.1cm}
            \State \ \ \ \ $w_{s}^{k+1} \leftarrow w_{s}^{k} + \rho(\hat{x}^{k}_s - \bar{x}^{k})$\vspace*{0.15cm}   
    \EndWhile
    \State \textbf{Return} solution with the lowest UB so far\vspace*{0.0cm}
\end{algorithmic}
\end{algorithm}
%
\section{Case study - The Brazilian power system}\label{sec:results}
This section will study the proposed a-BDMM method to solve the GEP for the Brazilian power system. The Brazilian power system is interconnected by a transmission network comprising 50 transmission lines connecting subsystems. The system has 550 thermal plants, 150 renewable plants, 200 hydro plants, 10 batteries, and 35 buses. We considered monthly stages (time steps), with three typical days per stage (week, weekend, and critical days), {\color{myblue}each of which comprised of 24 hours}. For each hydro unit, an affine policy with 24 coefficients is considered. The database configuration is based on a joint work of the system planning company (representing the Brazilian Energy Ministry) and private companies (market players and consulting companies) to study the impact of the large integration of renewables in the Brazilian system \cite{gizstudy}. The main assumptions are:

\begin{itemize}
    \item The PDE 2026 final system configuration is used as a starting point (see in \cite{EPE2020}).
    \item A target year is used where the demand is considered to be twice the demand of 2017, which amounts to 166 GW/1200 TWh (peak/annual energy).
    \item Annualized investment costs for the target year are used.
\end{itemize}
The complete data set used in this work can be found at \cite{Dados}. The complexity of the co-optimization of discrete investment decisions and long-term hydrothermal operational actions (optimal reservoir management) requires improved methods to enable assessments with large details from both long- and short-term uncertainties. Therefore, we compare the proposed a-BDMM method against the TBD, the BDMM, and the deterministic-equivalent (DE) formulation \eqref{eq:bigstart}-\eqref{eq:bigend} solved directly through a MILP algorithm. We compare the algorithm in terms of number of iterations and computational time. The analyses include 10 instances considering 5, 10, ..., 50 scenarios in $\Omega$ representing uncertainties in renewable generation, hydro inflows, and demand. Furthermore, we also analyze the effects of the nonanticipativity constraints in the expansion planning.

We used the following relevant parameters for the algorithm: GAP tolerance of 0.1\%; the maximum number of Benders iterations equal to 200; 12 stages representing months; and three typical days (weekdays, weekend days, and a critical day) for each stage. Because the proposed method significantly benefits from parallelism, we used the same computational resources to compare all methods. We used the Xpress solver (FICO, optimizer version 34.01) and an Amazon EC2 c5.12xlarge computer (48 processors and 96 GB of RAM).

It is relevant to highlight that vector $x$ comprises very different components, with different images and different weights in the original objective function \eqref{eq:bigstart}. While investment decisions are binary and already appear in the original objective function weighted by investment costs, LDR coefficients are real numbers and do not participate in the original objective function. Therefore, in our implementation, we modified the penalization term $\frac{\rho}{2}\|x-\bar{x}^{k}\|^2$ from expression \eqref{eq:master2start} to consider different penalty weights for components in $x^{INV}$ and $(x^{LDR}_{t,0}, x^{LDR}_{t})$. Therefore, we used the following quadratic penalty term:

{\color{myblue}
\begin{align}
&\frac{1}{2}\|x^{INV} - \overline{x}^{INV}\|^2_{diag(\mathcal{I})} + \notag \\
&\hspace{0.5cm} \frac{1}{2}\|x_{t}^{LDR} - \overline{x}_{t}^{LDR}\|^2 + \frac{1}{2}\|x_{t,0}^{LDR} - \overline{x}_{t,0}^{LDR}\|^2, \label{ModifiedPenalty}
\end{align}}
{\color{myblue}where, $||x||_P^2= {x^{\top} P x}$, when $P$ is the identity we omit the subscript. $diag(\cdot)$ converts a vector into a diagonal matrix. Hence, the first term is weighted by investment costs $\mathcal{I}$}. In this context, the quadratic deviation of investment decisions is penalized with half of their original weight (investment costs), whereas deviations of LDR coefficients are penalized with $\rho = 1$. This selection strategy constitutes a selection rule that improved the algorithm's efficiency (in terms of iterations) for all tested instances.

In the following sections, we analyze the performance, in terms of iterations and computational time, of the proposed method for different instance sizes, and the impact of nonanticipativity in the investment decisions, total cost and spot prices.
\subsection{Analysis of the decomposition algorithm}\label{sec:brazil-decomp}
Table \ref{tab:macro-results} summarizes some macro results of the proposed a-BDMM method for some selected instances. The number of constraints and variables are referring to the size of the DE version of the problem (instance \eqref{eq:bigstart}-\eqref{eq:bigend}). 

\begin{table}[h]
\caption{Results of the proposed algorithm for some instances (identified by the number of scenarios considered)}
\label{tab:macro-results}
\small
\begin{tabular}{lccccc}
\hline
\textbf{$\big| \Omega \big|$}& \textbf{15}    & \textbf{25}    & \textbf{30}      & \textbf{50}     \\ \hline
Constraints ($10^6$)                       & $17.5$ & $29.2$ & $35.0$ & $58.4$ \\
Variables  ($10^6$)                        & $24.5$ & $40.8$ & $49.2$ & $81.6$ \\
Execution time (min)                          & 28                & 65               & 120                & 206                \\
Number of iterations                        & 36                 & 33                 & 37                 & 31                 \\
Upper bound (M\$)                           & 71,423             & 57,087             & 42,600             & 30,036             \\
Lower bound (M\$)                           & 71,359             & 57,053             & 42,562             & 30,012             \\
Optimality GAP (\%)                         & 0.09              & 0.06              & 0.09               & 0.08               \\ \hline
\end{tabular}\\[.03cm]
\end{table}

Table \ref{tab:dattab-br} shows, for each instance and method, the number of iterations required to achieve a GAP of 0.1\%. Note that the number of iterations required by the proposed BDMM method is always smaller than that required by the TBD. Furthermore, the a-BDMM method further reduced this number, showing that the PH consensus scheme is effective in reducing the number of Benders loops needed to achieve the required GAP. Indeed, the a-BDMM required, on average (over all instances), 53\% fewer iterations than the BDMM and 68\% fewer iterations than the TBD. So, results show that the proposed a-BDMM method outperforms the benchmarks in terms of Benders iterations needed to solve the GEP for the Brazilian power system.

Table \ref{tab:dattab-br} also shows that the DE converges faster than the TBD and the proposed BDMM as long as the computer's memory is enough to address the problem. However, for larger instances in which MILP solvers fail to address the DE, Benders' approaches are still capable of providing high-quality solutions. Notwithstanding, it is important to highlight that the proposed a-BDMM breaks this pattern, achieving the required GAP faster than all methods for the larger eight out of ten tested instances. After considering 25 scenarios, the DE method fails to load the problem. Additionally, the same pattern observed for the number of iterations was observed for the computational time. The BDMM outperformed the TBD, and the a-BDMM outperformed the BDMM. Indeed, the a-BDMM is, on average (over all instances), 46\% faster than the BDMM and 60\% faster than the TBD. {\color{myblue}PH did not converge to the target gap after 500 minutes, thus we omitted PH from table III.}

\begin{table}[h]
\footnotesize
\centering
\caption{Number of iterations and computational time for different instance sizes and methods}
\label{tab:dattab-br}
\begin{tabular}{cccc|cccc}
\hline
     \multicolumn{4}{c}{\textbf{\scriptsize Number of Benders loops (iterations)}} & \multicolumn{4}{c}{\textbf{\scriptsize Execution time (in minutes)}} \\ \hline
    \textbf{$\big| \Omega \big|$}& \textbf{\scriptsize TBD} & \textbf{\scriptsize BDMM} & \textbf{\scriptsize a-BDMM} & \textbf{\scriptsize TBD} & \textbf{\scriptsize BDMM} & \textbf{\scriptsize a-BDMM} & \textbf{\scriptsize DE}   \\ \hline
5   & 176      & 78           & 34            & 29          & 17           & 11            & 7             \\       
10  & 132      & 82           & 33            & 48          & 33           & 17            & 11            \\       
15  & 136      & 89           & 36            & 88          & 59           & 28            & 31            \\       
20  & 125      & 83           & 36            & 132         & 88           & 46            & 59            \\       
25  & 113      & 83           & 33            & 193         & 142          & 65            & -             \\       
30  & 101      & 72           & 37            & 276         & 203          & 120           & -             \\       
35  & 87      & 64           & 34            & 287         & 226          & 135           & -             \\       
40  & 73      & 61           & 32            & 252         & 237          & 131           & -             \\       
45  & 88      & 70           & 33            & 356         & 320          & 163           & -             \\       
50  & 78      & 63           & 31            & 443         & 407          & 206           & -             \\ \hline
\end{tabular}
\end{table}

For comparison purposes, the algorithm a-BDMM presented a total execution time of 206 minutes for the 50-scenario instance, where 21\% (44 minutes) was used to solve the master problems and 79\% (162 minutes) to solve the subproblems.

\begin{figure}[h]
\includegraphics[scale=0.4]{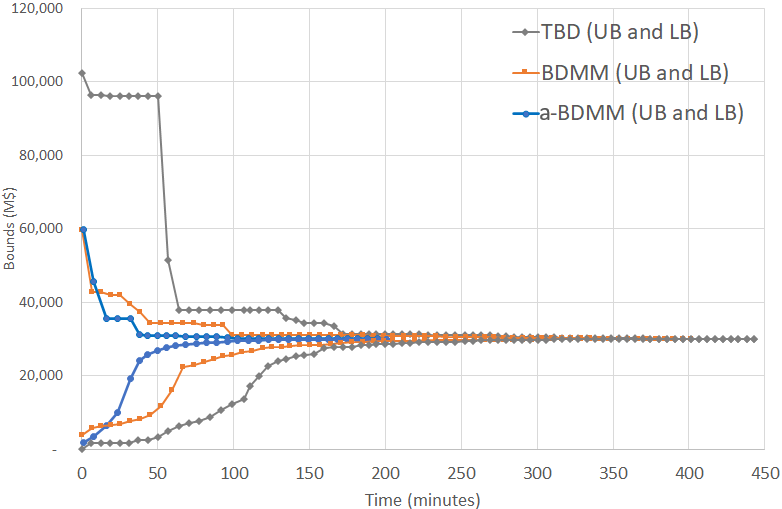}
\caption{Lower and upper bound for each algorithm}
\label{fig:convergence}
\end{figure}

\begin{figure}[h]
\centering
\includegraphics[scale=0.4]{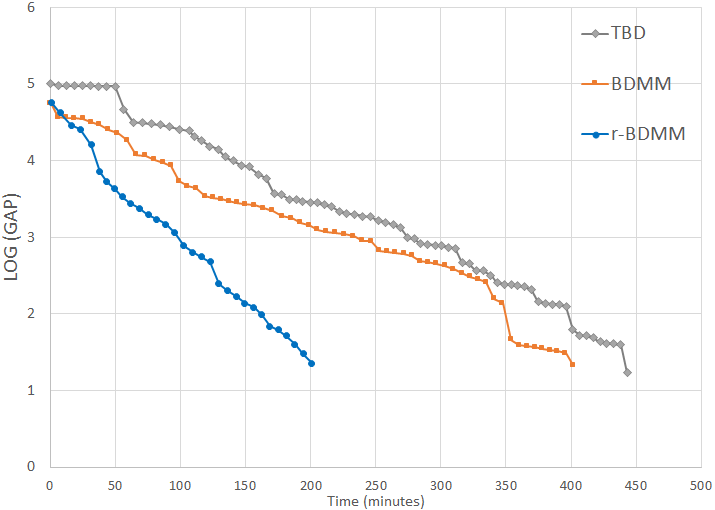}
\caption{{\color{myblue}Log of absolute GAP for each algorithm}}
\label{fig:convergence_log}
\end{figure}

Finally, Figure \ref{fig:convergence} and {\color{myblue}Figure \ref{fig:convergence_log}} compare the convergence over time of the a-BDMM, BDMM, and TBD, for the instance with 50 scenarios. It's clear that the convergence of the a-BDMM outperforms the BDMM and TBD. Also, after 1 hour of running time, the GAP values are: 9.76\% for the a-BDMM; 53.15\% for the BDMM; and 490\% for the TBD. The GAP after 2 hours decreases to: 1.59\% for the a-BDMM; 11.75\% for the BDMM; and 68.3\% for the TBD. These results corroborate the superiority of our proposed method to solve the GEP for the Brazilian system.

\vspace{-0.2cm}
\vspace{-0.2cm}
\subsection{Benefits of a nonanticipative operational policy}\label{sec:ldr}

This section analyzes the benefits of considering a multistage (nonanticipative) operational policy when deciding the investment plans for the Brazilian power system. To do that, we consider two cases:

\vspace{0.1 cm}
\noindent \textbf{Nonanticipative policy} -- We solve problem \eqref{eq:bigstart}-\eqref{eq:bigend} with the a-BDMM algorithm and $|\Omega_{50}|=50$ scenarios (the same case study analyzed in Section \ref{sec:brazil-decomp}). The operational results obtained in the optimization will be referred to as \textit{in-sample}. Then, we fix the optimal value obtained for $x^*$ and evaluate $Q(x^*)$ with $|\Omega_{1000}|=1000$ scenarios (i.e., we evaluate the operational part of the problem, \eqref{eq:recourse}, with one thousand unseen scenarios). These results will be referred to as "out-of-sample".
    
\noindent \textbf{Anticipative policy} -- We solve the same problem \eqref{eq:bigstart}-\eqref{eq:bigend}, disregarding the nonanticipative constraints \eqref{mod:affine} for $|\Omega|=50$. In this context, we are considering an anticipative approximation for the hydrothermal dispatch costs when deciding the generation investment plans. Then, as in the previous case (Nonanticipative), we will provide results for both in-sample and out-of-sample cases. However, in order to assess the benefit of a multistage nonanticipative operational policy when making investment decisions in hydrothermal power systems, the out-of-sample analysis must be carried out based on an implementable (nonanticipative) policy. To do that, we fix the optimal investment part of the solution found with the anticipative approximation, i.e., $x^{INV*}_A$, load the nonanticipative constraints to the in-sample problem, and solve it again to define $x^{LDR*}(x^{INV*}_A)$. Then, with the complete first-stage vector, $x^*_A = [x^{INV*}_A \; x^{LDR*}(x^{INV*}_A)]$, we evaluate the out-of-sample operational cost. To do that, we use the same 1000 scenarios used in the out-of-sample evaluation of the Nonanticipative policy.

Notwithstanding, it is clear that the Anticipative policy is motivated by its lower computational burden. In this case study, the in-sample optimization of the Nonanticipative policy took 206 min, whereas the Anticipative case took 68 minutes. Therefore, this work aims to highlight the benefits of considering a nonanticipative operational policy to justify its higher computational times. In our comparison, we used the following metrics:

\begin{itemize}
    \item \emph{Expected total cost} -- assessed with out-of-sample nonanticipative operational results.
    \item \emph{Regret} -- measured as the difference between the in-sample and out-of-sample expected total cost.
    \item \emph{The average spot price} -- assessed with out-of-sample operational results.
    \item \emph{The 95th-percentile of the spot price} -- assessed with out-of-sample operational results.
    \item \emph{The average uncertainty level of spot price} -- assessed with out-of-sample operational results. This metric is defined as (see \cite{soares2017solution}): $\dfrac{1}{|\mathcal{T}|}\sum\limits_{t=1}^{|\mathcal{T}|}\left(\mathbf{Q}_t^{95\%}-\mathbf{Q}_t^{5\%}\right)$, where $\mathbf{Q}_t^{\alpha\%}$ represents the $\alpha\%$ quantile of a given variable at stage $t$.
    \item \emph{The time variability of the spot price} -- assessed with out-of-sample operational results. This metric is defined as (see \cite{soares2017solution}):  $\dfrac{1}{|\mathcal{T}|-1}\sum\limits_{\omega\in\Omega_{1000}}p_\omega\sum\limits_{t=2}^{|\mathcal{T}|}\left| \dfrac{\pi_{t,\omega}-\pi_{t-1,\omega}}{\pi_{t-1,\omega}}\right|$.
    \item \emph{Value of the nonanticipative policy} -- the difference between the expected total cost, assessed with out-of-sample operational results, of the Anticipative and Nonanticipative policies, i.e., $VNAP = \mathcal{I}^{\top}x^*_A + Q_{\Omega_{1000}}(x^*_A) - (\mathcal{I}^{\top}x^* + Q_{\Omega_{1000}}(x^*))$.   
\end{itemize}

Table \ref{tab:affine} shows the in-sample and out-of-sample costs for both Anticipative and Nonanticipative policies. Table \ref{tab:affine} shows that the anticipative policy, albeit 10,364M\$ cheaper than the Nonanticipative when analyzed with in-sample results, is actually 6,579M\$ (or 8.27\%) more expensive when analyzed with out-of-sample results. This difference defines the benefit or value of considering a Nonanticipative policy when making the investment decisions: $VNAP = 86,116 - 79,537 = 6,579M\$$ in absolute terms, or 7.64\% of the total cost and $16.18\%$ of the investment cost obtained with the Anticipative policy. The difference between what was expected when optimizing and what we got when actually implementing the solutions defines the regret metric, which values 20,679 M\$ (or 24\% of the total cost and 50.85\% of the investment cost obtained with the anticipative operational policy). 

\begin{table}[h]
\centering
\caption{Impact of the stochastic policy in terms of total cost}
\label{tab:affine}
\small
\begin{tabular}{cccc}
\hline
  \multicolumn{4}{c}{\textbf{In-sample results}} \\ \hline
  \textbf{Policy} & \thead{\textbf{Investment cost} \\ \textbf{(M\$)}}& \thead{\textbf{Operational cost} \\ \textbf{(M\$)}} & \thead{\textbf{Total cost} \\ \textbf{(M\$)}} \\ \hline
Anticipative                             & 40,662          & 24,775          & 65,437 \\
Nonanticipative                             & 42,646          & 33,155          & 75,801 \\
Difference                           & -1,984          & -8,380         & -10,364  \\ \hline
  \multicolumn{4}{c}{\textbf{Out-of-sample results}}\\ \hline
  \textbf{Policy} & \thead{\textbf{Investment cost} \\ \textbf{(M\$)}}& \thead{\textbf{Operational cost} \\ \textbf{(M\$)}} & \thead{\textbf{Total cost} \\ \textbf{(M\$)}} \\  \hline
Anticipative                             & 40,662          & 45,454          & 86,116 \\
Nonanticipative                             & 42,646          & 36,891          & 79,537 \\
Difference                           & -1,984          & 8,563          & 6,579  \\ \hline

\end{tabular}
\end{table}

\begin{table*}[h]
\centering
\caption{Out-of-sample metrics for the temporal inconsistency}
\label{tab:metrics}
\small
\resizebox{\textwidth}{!}{%
\begin{tabular}{ccccccc}
\hline
  \textbf{Policy} & \thead{\textbf{Total costs} \\ \textbf{(M\$)}} & \thead{\textbf{Regret} \\ \textbf{(M\$)}} & \thead{\textbf{Average spot} \\ \textbf{price (R\$/MWh)}} & \thead{\textbf{95th-percentile of} \\ \textbf{the spot price (R\$/MWh)}} & \thead{\textbf{Average uncertainty level} \\ \textbf{of the spot price (R\$/MWh)}} & \thead{\textbf{Time variability} \\ \textbf{of the spot price (\%)}} \\ \hline
Anticipative                             & 86,116          & 20,679          & 335.4 & 900.9 & 1276.7 & 105\\
Non Anticipative                             & 79,537          & 3,736          & 156.2 & 542.9 & 478.0  & 65\\
Difference                           & 6,579          &       16,943    & 179.2& 358.0 & 798.7 & 40  \\ \hline

\end{tabular}
}
\end{table*}

\begin{figure}[h]
\centering
\includegraphics[scale=0.3]{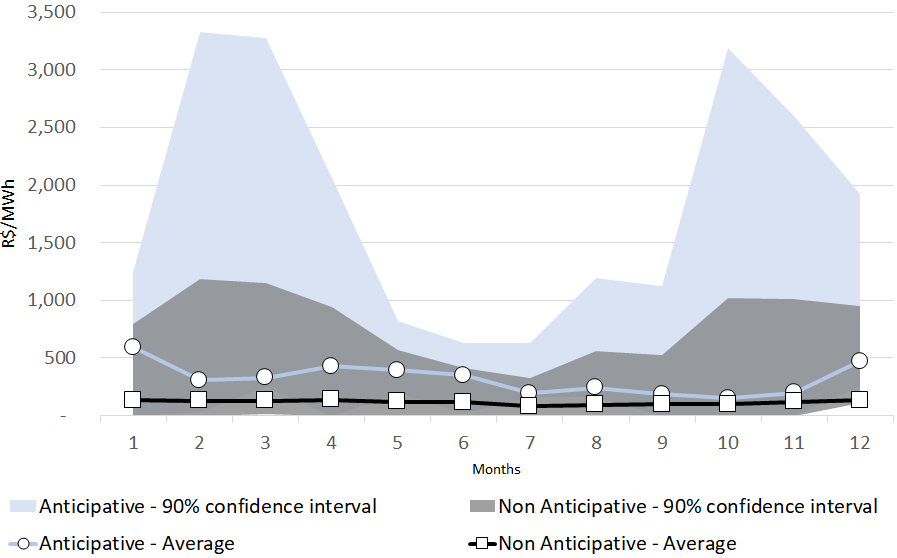}
\caption{90\% confidence interval of energy spot prices for both Antecipative and Nonanticipative policies}
\label{fig:price}
\end{figure}

\begin{figure}[h]
\centering
\includegraphics[scale=0.35]{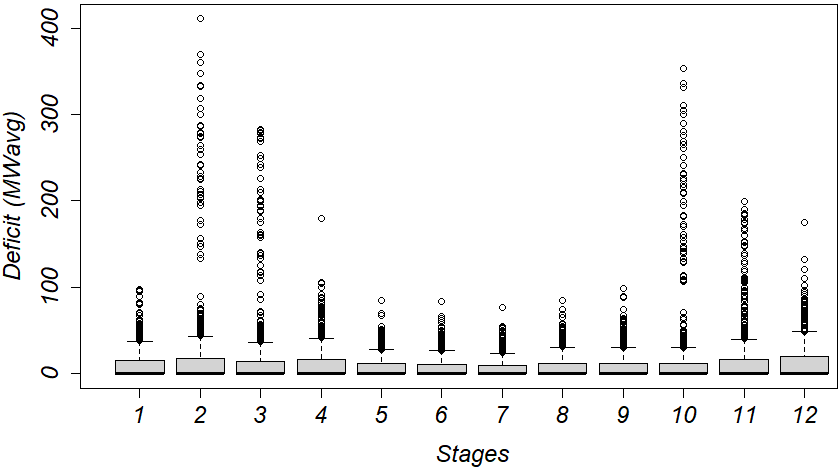}
\caption{Variability and amount of deficit for the Anticipative policy}
\label{fig:boxplot}
\end{figure}

Figure \ref{fig:price} shows the 90\% confidence interval for the spot price for both policies. Since both cases considered interconnections between areas, the spot prices for each of the regions are exactly the same. This figure shows that the investments made under an anticipative policy, when actually operating the system, produce much higher and uncertain spot prices. Furthermore, Table \ref{tab:metrics} shows the metrics presented at the beginning of this section. We can see that, besides being 8.27\% more expensive, the non-implementable policy brings higher spot price on average and for high quantiles, higher volatility (uncertainty level) and higher temporal variability. Additionally, the deficit risk for the Nonanticipative policy achieved a 0\% probability in the out-of-sample, while the Anticipative policy exhibited several scenarios with deficit as shown in Figure \ref{fig:boxplot}.

These results are consistent with the results obtained in previously reported works where simplifications were used in the opportunity cost assessment \cite{street2020assessing,brigatto2017assessing}. In this work, however, we extend this idea to the expansion planning level. 

\vspace{-0.2cm}
\section{Conclusions}
\label{sec:conclusions}
This work presented a novel Benders Decomposition with multiple master (BDMM) problems. We show that the proposed method significantly improves the performance of the traditional Benders decomposition by using different parallel master problems, each of which considering a primal cut associated with a given scenario. Furthermore, we also propose a novel acceleration approach for our BDMM (a-BDMM). The proposed approach is based on the consensus idea behind the Progressive Hedging method and is capable of significantly reducing convergence times when applied to the hydrothermal generation expansion planning problem.

In this work, we studied the Brazilian hydrothermal power system, which highly relies on the assessment of the opportunity cost of water through multistage nonanticipative operational policies. We show that the consideration of a multistage nonanticipative policy, rather than the less computationally intensive anticipative approximation, brings relevant benefits to the optimal investment decisions. The proposed a-BDMM has shown to be crucial for solving the large-scale Brazilian case.

Based on realistic data from the Brazilian power system (see \cite{Dados} based on \cite{gizstudy}), the case studies presented in this work allow us to convey the following concluding remarks:

\begin{itemize}
    \item The proposed BDMM method outperforms the traditional Benders decomposition benchmark in both the number of iterations (30\% on average) and computational time (23\% on average) under the same computational resources.
    \item The proposed a-BDMM method outperforms the BDMM in both number of iterations (53\% on average) and computational time (46\% on average). So, the proposed a-BDMM outperforms the traditional Benders decomposition in terms of the number of iterations (68\% on average) and computational time (60\% on average).
    \item The proposed a-BDMM can solve large-scale instances of the Brazilian expansion planning problem considering the co-optimization of discrete investment decisions and nonanticipative (multistage stochastic) long-term hydrothermal operational actions to manage reservoir levels.
    \item The value of considering the nonanticipative hydrothermal operational policy (multistage dispatch under uncertainty) when deciding the investment plans is 6,579 M\$, i.e., 16.18\% of the investment costs obtained with the Antecipative policy.
    \item Regarding the spot price profile, the consideration of a nonanticipative operational policy brings other benefits in comparison to the anticipative counterpart, namely, (i) a reduction of 53.4\% on the expected annual prices, (ii) a reduction of 39.7\% on the 95th-percentile of the annual  prices, (iii) a reduction of 62.6\% on the price uncertainty metric; and (iv) a reduction of 38\% on the temporal variability metric.
    \item The regret of considering an anticipative (non-implementable) operational policy in the expansion problem is 20,679 M\$, which is 51\% of the investment cost.
\end{itemize}
\ignore{\cite{koltsaklis2015multi} \cite{pina2013high} \cite{li2018robust} \cite{liu2018multistage}  \cite{thome2019stochastic}}
\vspace{-0.35cm}
\bibliographystyle{IEEEtran}
\bibliography{IEEEabrv,Bibliography}

\end{document}